\documentclass[12pt]{article}
\usepackage[russian,english]{babel}
\usepackage[T2A]{fontenc}           
\usepackage[cp1251]{inputenc}       
\usepackage{amsfonts}
\usepackage[active]{srcltx}
\usepackage{amsmath}
\usepackage{amssymb}
\usepackage{cite}
\usepackage{graphicx}
\usepackage{color}
\newenvironment{sciabstract}{%
\begin{quote} \footnotesize} 
{\end{quote}}
\newtheorem{theorem}{Theorem}

\newtheorem{corollary}{Corollary}

\newtheorem{example}{Example}
\newtheorem{lemma}{Lemma}

\newtheorem{remark}{Remark}
\usepackage{graphicx}
\title{\bf Solving ill-conditioned linear algebraic systems using methods that improve conditioning
}
\author
{A.S. Leonov$^1$\\
\\
\normalsize{$^{1}$National Research Nuclear University MEPhI,}\\
\normalsize{Moscow, Russian Federation}\\
\normalsize{E-mail: asleonov@mephi.ru}
}
\date{}

\begin{document}

\maketitle

\begin{sciabstract}
We consider the solution of systems of linear algebraic equations (SLAEs) with an ill-conditioned or degenerate exact matrix and an approximate right-hand side. An approach to solving such a problem is proposed and justified, which makes it possible to improve the conditionality of the SLAE matrix and, as a result, obtain an approximate solution that is stable to perturbations of the right hand side with higher accuracy than using other methods. The approach is implemented by an algorithm that uses so-called minimal pseudoinverse matrices. The results of numerical experiments are presented that confirm the theoretical provisions of the article.


\end{sciabstract}


\section{Introduction}
When studying various scientific problems, it is often necessary to solve systems of linear algebraic equations (SLAEs) of very large dimensions. Such systems often turn out to be ill-conditioned or degenerate. As a consequence, classical methods of linear algebra for solving these SLAEs (see, for example, \cite{fox,fors,lowson,golub1}) sometimes turn out to be of little use due to instability when rounding errors accumulate. This problem can be partially solved by using calculations on computers with double, quadruple and higher precision. However, this approach is not applicable to all computing systems on personal computers (PCs). Therefore, the practical solution of ill-conditioned and degenerate SLAEs requires the use of special methods customized to combat numerical instability.

Currently, many such methods are known with varying complexity and scope of applicability (see \cite{lowson,voev,tyrt}, etc.). Let us especially mention here methods that are relevant to the subsequent presentation. Among them, the regularization method of A.N. Tikhonov and its variant in the form of the regularized least squares method \cite{tikh1} are widely known. The TSVD (truncated singular value decomposition) method in its various modifications is also often used (see, for example, \cite{TSVD1,TSVD2,golub2}). A minimal pseudoinverse matrix (MPM) method has also been developed (see \cite{leonov1} and \cite{tikh2,leonov_book}) that will be mainly used further.

A distinctive feature of this work is the construction of special methods for solving SLAEs, \emph{significantly improving the condition number} of their matrix (i.e., significantly reducing them). This is how the proposed methods differ from those listed above, in which reducing the condition number is not the goal, but turns out to be a kind of by-product.

The article is organized as follows. In Sect.2, methods for stable solution of SLAEs, which are used to varying degrees in the work, are briefly described. In Sect.3, some auxiliary statements are considered. The problem statement, the proposed solution scheme and its justification are presented in Sect.4. A specific method that implements the proposed scheme for solving SLAEs is proposed in Sect.5. The results of numerical experiments with model problems and conclusions are presented in Sect.6.

\section{On some methods for stable solution of SLAEs} 
Let $\mathbb{R}^{n}_c, \mathbb{R}^{m}_c$ be spaces of column vectors with Euclidean norms. We will solve SLAE of the form
\begin{equation}\label{SLAU1}
Az=u,\,\ z\in \mathbb{R}^{n}_c,\,u\in \mathbb{R}^{m}_c
\end{equation}
with a real matrix $A,\,\|A\|\neq 0$, $\dim A=m\times n$, and right-hand side $u\neq 0$. Since this system may not be solvable in
in the classical sense, we will consider its pseudosolutions $z^{\ast}\in \mathbb{R}^{n}_c$, i.e. elements satisfying equality
\begin{equation*}
\left\| Az^{\ast }-u\right\| =\inf \{\left\| Az-u\right\| :\,\
z\in \mathbb{R }^{n}_c\}.
\end{equation*}
We denote the set of all pseudo-solutions for the exact data of the problem $\{A,u\}$ as $Z^{\ast }$. For a compatible system, these pseudo-solutions represent conventional solutions. For any system of the form \eqref{SLAU1} there is a unique normal pseudosolution $\bar{z}$, i.e. such a pseudo-solution for which
\begin{equation*}
\left\| \bar{z}\right\| =\inf \{\left\| z\right\| :\,z\in Z^{\ast
}\}.
\end{equation*}
For a solvable system \eqref{SLAU1} it is called a normal solution. For a uniquely solvable system, the normal solution coincides with the classical solution. Our goal is to find the vector $\bar z$ from the data $\{A,u\}$. The formal solution to this problem is given using the pseudo-inverse matrix $A^+$: $\bar z = {A^ + }u$ (see, for example, \cite{alb,lowson,golub1}).

Often this formulation of the problem is modified. It is believed that instead of exact data $\{A,u\}$ we have at our disposal their approximations $\{A_{h},u_{\delta }\}$, $\mathrm{dim} A_h=\mathrm{dim } A$, $u_{\delta }\in \mathbb{R}^{m}_c$, and error levels, i.e. numbers $\eta =(h,\delta )$, such that $\left\| A_{h}-A\right\| \leq h,\,\left\| u_{\delta }-u\right\| \leq \delta $. Here and in what follows, the norms of matrices are considered Euclidean: $\left\| A\right\| = {\left\| A\right\|_E}$. It is required to construct approximations $z_{\eta }\in \mathbb{R}^{n}_c$ to the solution $\bar{z}$, which have the stability property: $\left\| z_{\eta}-\bar{z} \right\| \rightarrow 0$ at $\eta \rightarrow 0$.

The solution to this problem using the formula ${z_\eta } = A_h^+ {u_\delta }$ turns out to be unstable in the general case, because convergence of $A_h^+ \to {A^+ }$ for $h\to 0$ is not guaranteed. Therefore, other methods have been proposed. For example, we can take as a solution the element $z_\eta=z_\eta ^{\alpha (\eta )}$, which is obtained by A.N. Tikhonov’s regularization method:
\[z_\eta ^{\alpha (\eta )} = {(\alpha (\eta )I + A_h^*{A_h})^{ - 1}}A_h^*{u_\delta }=T_1(A_h)u_\delta.\] Here the parameter $\alpha (\eta )>0$ is chosen by algorithms described in \cite{tikh1,tikh2,leonov_book,tikh3,moroz} and other works. These algorithms provide the required convergence of approximations. To implement the regularization method, it is necessary to know the errors levels of the matrix and the right-hand side of the SLAE.

Another solution to the problem can be obtained using the TSVD method. It is based on the artificial zeroing of small singular values of the system matrix. The method implementation scheme is as follows. The singular value decomposition (SVD) of the approximate matrix is found:
\[A_{h}=U_hR_hV^{T}_h,\,R_h=\mathrm{diag}(\rho _{1},\rho _{2},...,\rho
_{M}),\,\rho _{1}\geq \rho _{2}\geq ...\geq \rho _{M}\geq 0\] with
$M=\min (m,n)$. Next, using some method, the number
$\kappa\in \mathbb{N}$ is selected that plays the role of a regularization parameter for which $\rho _{\kappa}>0$ and  $\rho _{k}=0$ for $k=\kappa+1,\kappa+2,...,M$. Then a new matrix of singular numbers
$R_{\kappa}=\mathrm{diag}(\rho _{1},...,\rho _{\kappa},0,...,0)$  is constructed. From it, as well as from the orthogonal matrices $U_h,V_h$ of the SVD, an approximate solution of the SLAE is found:
\begin{equation*} z_{\kappa }=V_hR_{\kappa }^{+}U^{T}_hu_{\delta }=T_2(A_h)u_{\delta },\,\
R_{\kappa }^{+}=\mathrm{diag}(\rho _{1}^{-1},...,\rho _{\kappa
}^{-1},0,...,0).
\end{equation*}
This corresponds to using, instead of the matrix $A_{h}$, another approximate matrix $\hat{A}_{h}=U_hR_{\kappa}V^{T}_h$ and its pseudo-inversion. The question of what is considered a small singular number, i.e. how to choose $\kappa$ is decided differently. For example, in the article \cite{golub2} it is proposed to choose $\kappa=\kappa(h)$ so that
\begin{equation*}
(\rho _{\kappa +1}^{2}+\rho _{\kappa +2}^{2}+...+\rho
_{M}^{2})^{1/2}\leq h,\,\ (\rho _{\kappa }^{2}+\rho _{\kappa
+1}^{2}+...+\rho _{M}^{2})^{1/2}>h.
\end{equation*}
This corresponds to choosing from all matrices $R_{\kappa }$ satisfying the condition $\left\|R_h-R_{\kappa }\right\|_E\leq h$ the matrix that has the minimum rank. This method requires knowledge of the matrix error $h$, but does not use the $\delta$ value.

Both considered methods are based on stable approximation of the pseudoinverse matrix $A^+$ by matrices $T_{1,2}(A_h)$. Note that the regularization method uses a full-rank matrix $T_{1}(A_h)$, while the matrix $T_{2}(A_h)$ of the TSVD method may have a significantly lower rank than the original matrix $A_h$. The condition number of matrices $T_{1,2}(A_h)$ is also important. In what follows we use the spectral condition numbers $\nu _{s}(A)=\rho _{1}/\rho _{r}$, where $r$ is the rank of the corresponding matrix ${A}$. We will compare the values $\nu _{s}(T_{1,2}(A_h))$ and the condition numbers of matrices arising in the methods proposed below.

In the works \cite{leonov1,tikh2,leonov_book}, a method of minimal pseudoinverse matrix (MPM) was proposed and studied, which not only allows one to stably calculate $A^+$, but also has a number of optimal properties. Let us briefly describe this method. We introduce the space $\mathcal{A}$, consisting of matrices of size $m\times n$ with Euclidean norms. Let $\bar A\in \mathcal{A}$ be an unknown exact matrix for which we need to find a pseudo-inverse. Approximate data for this are
known quantities $A_h\in \mathcal{A}$ and $h$ such that $\left\|A_h-\bar A\right\|\leq h$. Next, we define a set of matrices
$\mathcal{A}_{h}=\{A\in \mathcal{A}:\,\left\| A_{h}-A\right\| \leq
h\}$, comparable to the matrix $A_h$ in accuracy. It is clear that $\bar A\in \mathcal{A}_h$. Consider the following extremal problem: find a matrix $\tilde{A}_{h}\in \mathcal{A}$ for which equality
\begin{equation}\label{SLAU200}
\left\| \tilde{A}_{h}^{+}\right\| =\inf \{\left\| A^{+}\right\|
:\,A\in \mathcal{A}_{h}\}
\end{equation}
is satisfied.
The solution to the problem (\ref{SLAU200}) is called the \emph{matrix of the MPM method}, and the corresponding pseudoinverse matrix $\tilde{A}_{h}^{+}$ is a \emph{minimal pseudoinverse matrix}. It is proved that these matrices exist, and for $h\rightarrow 0$ the accuracy estimate is valid: $\left\| \tilde{A}_{h}^{+}-\bar{A}^{+}\right\|
\thicksim 2h\left\| \bar{A}^{+}\right\| ^{2}$. The estimate is asymptotically unimprovable. This means that the MPM method gives the optimal order of accuracy of approximation to the pseudoinverse matrix $\bar{A}^{+}$.

The MPM method is easily implemented numerically. The corresponding algorithm consists of the following steps (see \cite{leonov1},
\cite[Chapter 5]{tikh2}).

\noindent 1. Finding the SVD of given matrix
$A_{h}$: $A_{h}=U_{h}R_{h}V_{h}^{T}$, where $R_{h}=\mathrm{diag}(\rho _{1},\rho _{2},...,\rho _{M})$, $\rho _{1}\geq \ rho _{2}\geq ...\geq \rho _{M}\geq 0$, $M=\min (m,n)$, and $U_{h},V_{h}$ are orthogonal matrices of sizes $m\times m$ and $n\times n$.

\noindent 2. A) Finding the quantities $\lambda _{k}=\frac{27}{16}\rho _{k}^{4}$ and calculating the functions
\begin{equation*}
\rho _{k}(\lambda )=\{\rho _{k}x_{k}(\lambda ),\, 0<\lambda \leq
\lambda _{k};\,0,\, \lambda \geq \lambda _{k}\};\, \rho
_{k}(0)=\rho _{k},\,\,k=1,2,...,M,
\end{equation*}
for $\lambda \geq 0$. Here $x_{k}(\lambda )$ is a solution to the equation $ x^{4}-x^{3}=\lambda \rho _{k}^{-4}$, belonging to the segment $\lbrack 1, \frac{3}{2}]$.

\noindent B) calculation of the function
\begin{equation*}
\beta (\lambda )=\sum_{k=1}^{M}[\rho _{k}(\lambda )-\rho
_{k}]^{2}, \,\lambda\geq 0.
\end{equation*}

\noindent 3. Calculation of the solution $\lambda (h)>0$ of the equation $\beta (\lambda )=h^2$. In practice, steps 2 A), B) are used repeatedly here.

\noindent 4. Finding a regularized singular value matrix
\[
\tilde{R}_{h}=\mathrm{diag}(\rho _{1}[\lambda (h)],\rho
_{2}[\lambda (h)],...,\rho _{M}[\lambda (h)])
\]
and its pseudo-inverse
\[
\tilde{R}_{h}^{+}=\mathrm{diag}\{\theta \lbrack \rho _{1}[\lambda
(h)]],\theta \lbrack \rho _{2}[\lambda (h)]],...,\theta \lbrack
\rho _{M}[\lambda (h)]]\},
\]
where $\theta (\rho )=\{\rho ^{-1},\, \rho >0;\,0,\, \rho =0\}$.

\noindent 5. The final calculation of the minimum pseudoinverse matrix: $\tilde{A}_{h}^{+}=V_{h}\tilde{R}_{h}^{+}U_{h}^{T}$ and,
if necessary, finding the matrix of the MPM method
$\tilde{A}_{h}=$ $ U_{h}\tilde{R}_{h}V_{h}^{T}$.

As a commentary on the algorithm, we note that the function $\beta (\lambda )$ increases monotonically, having as its breakpoints the numbers $\lambda _{k}$. At these points the uniqueness of the function $\beta (\lambda )$ is violated, and at them it takes the values $\beta (\lambda _{k}-0)$ and $\beta (\lambda _{k}+0)$.
Therefore, when speaking about solving the equation $\beta (\lambda )=h^2$, we mean finding a generalized solution, i.e. or jump points of the function $\beta (\lambda )$ through the value $h^2$, or an ordinary root. The issues of implementing the MPM algorithm, as well as the influence of errors in calculating singular numbers for the original matrix on the algorithm, are discussed in more detail in \cite{leonov1}, \cite{tikh2}. It should also be noted that the MPM method produces a matrix $\tilde{A}_{h}$ that is optimal in order of the condition number.The method does not require knowledge of $\delta$. Below we will use minimal pseudoinverse matrices in a new algorithm for solving SLAEs, which guarantees an improvement in the condition number of the system matrix.

\section{Auxiliary statements}
In this section we assume that for the exact matrix $\bar{A}\in \mathcal{A}$, $\|\bar{A}\|\neq 0$, a family of approximate matrices $\tilde{A}_h\in \mathcal{A}$, $\forall\, h\ge 0 $ is given, such that the approximation condition $\left\Vert \tilde{A}_{h}-\bar{A}\right\Vert \leq h$ is satisfied.
\begin{lemma} \label{lem1}
Let the matrices $\tilde{A}_h$ satisfy the relation $\left\Vert \tilde{A}_{h}^{+}\right\Vert \leq \left\Vert \bar{A}^{+}\right\Vert$, $\forall h,~0\leq h<\left\Vert \bar{A}^{+}\right\Vert ^{-1}$. Then:
A) $\mathrm{Rg}\tilde{A}_{h}=\mathrm{Rg}\bar{A}$; B)
$\left\| \tilde{A}_{h}^{+}-\bar{A}^{+}\right\| \leq h\left\| \bar{A}%
^{+}\right\| ^{2}(1-h\left\| \bar{A}^{+}\right\| )^{-3}$.
\end{lemma}
Proof. We denote as $\bar{\rho}_{k},\tilde{\rho}_{k},~1\leq k\leq M$ the singular values of the matrices $\bar{A},\tilde{A }_{h}$. Then (see \cite{leonov1}) the relations
\begin{equation}\label{eq1}
\sum\limits_{k=1}^{M}(\tilde{\rho}_{k}-\bar{\rho}_{k})^{2}\leq \left\Vert \tilde{A}_{h}-\bar{A}\right\Vert ^{2}\leq h^{2},~\left\Vert \tilde{A}_{h}^{+}\right\Vert ^{2}=\sum\limits_{k=1}^{{r}_h}\tilde{\rho}_{k}^{-2},~\left\Vert \bar{A}^{+}\right\Vert ^{2}=\sum\limits_{k=1}^{r}\bar{\rho}_{k}^{-2},
\end{equation}
are satisfied, where $r_h=\mathrm{Rg}\tilde{A}_{h}, r=\mathrm{Rg}\bar{A}$. Let us prove A), i.e. equality $r_{h}=r$ for $0\leq h<\left\Vert \bar{A}^{+}\right\Vert ^{-1}$. We first assume that $r_{h}<r$. Then $\tilde{\rho}_{k}=0$ for $k>r$, and therefore, taking into account \eqref{eq1} we obtain
\[\bar{\rho}_{r}\leq \sum\limits_{k=1}^{r-1}(\tilde{\rho}_{k}-\bar{\rho}_{k})^{2}+
\sum\limits_{k=r}^{M}\bar{\rho}_{k}^{2}=\sum\limits_{k=1}^{M}(\tilde{\rho}_{k}
-\bar{\rho}_{k})^{2}\leq h^{2},~\bar{\rho}_{r}^{-2}\leq \left\Vert \bar{A}^{+}\right\Vert^2. \]
It follows that $h^{2}\geq \left\Vert \bar{A}^{+}\right\Vert ^{2}$, and this contradicts the condition on $h$.

If we assume that $r_{h}>r$, then $\bar{\rho}_{r_{h}}=0$, and similarly we obtain
\[
\tilde{\rho}_{r_{h}}^{2}=(\tilde{\rho}_{r_{h}}-\bar{\rho}_{r_{h}})^{2}\leq \sum\limits_{k=1}^{M}(\tilde{\rho}_{k}-\bar{\rho}_{k})^{2}\leq h^{2},~\tilde{\rho}_{r_{h}}^{-2}\leq \sum\limits_{k=1}^{r_{h}}\tilde{\rho}_{k}^{-2}=\left\Vert \tilde{A}_{h}^{+}\right\Vert ^{2}\leq \left\Vert \bar{A}^{+}\right\Vert ^{2}.
\]
This again implies the inequality $h^{2}\geq \left\Vert \bar{A}^{+}\right\Vert ^{2}$, which contradicts the conditions. As a result we get that $r=\mathrm{Rg}\bar{A}=r_{h}=\mathrm{Rg}\tilde{A}_{h}$.

The proof of part B) is based on the application of the estimate established in \cite{leonov1}:
\[
\left\Vert A^{+}-B^{+}\right\Vert \leq \left\Vert A-B\right\Vert \left\Vert A^{+}\right\Vert ^{2}/\left( 1-\left\Vert A-B\right\Vert \left\Vert A^{+}\right\Vert \right) ^{3}
\]
for matrices with $\dim A=\dim B,~\mathrm{Rg}A=\mathrm{Rg}B$. Setting $A=\bar{A},~B=\tilde{A}_{h}$ in this estimate and taking into account the proven part A) of the lemma, we obtain the required inequality. $\square$
\begin{corollary}\label{cor1}
If $\tilde{A}_{h}\rightarrow \bar{A}$ and $\left\Vert \tilde{A}_{h}^{+}\right\Vert \leq \left\Vert \bar {A}^{+}\right\Vert $ for $h\rightarrow 0$, then $\tilde{A}_{h}^{+}\rightarrow \bar{A}^{+}$.
\end{corollary}
This leads to next corollary.
\begin{corollary}\label{cor2}
Let $\delta\to 0$, and for some function $h(\delta)\to 0$ a matrix $A_{h(\delta)}$ is found such that
\[\left\|A_{h(\delta)}-A\right\|\rightarrow 0,\,\left\| {A_{h(\delta )}^+ } \right\| \le \left\| {{A^+ }} \right\|.\] Then $A_{h(\delta )}^+ \to {A^+ }$ for $\delta\to 0$ and therefore ${z_\delta } = A_{h(\delta )}^+ {u_\delta } \to \bar z$.
\end{corollary}

\section{Detailed formulation of the problem, proposed solution scheme and its justification}\label{sec4}
\begin{figure}[h]
  \centering
\includegraphics[width=80mm,height=70mm]{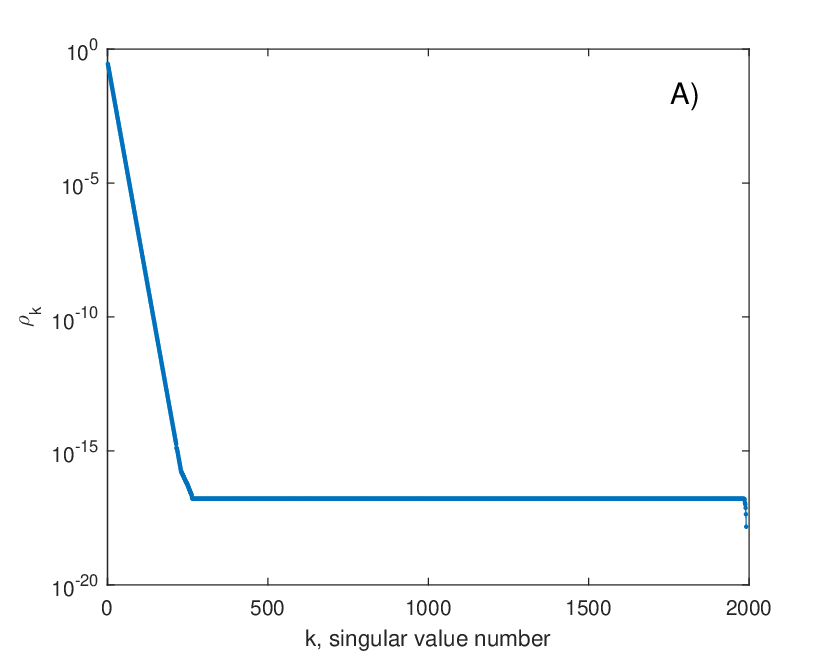}%
\includegraphics[width=80mm,height=70mm]{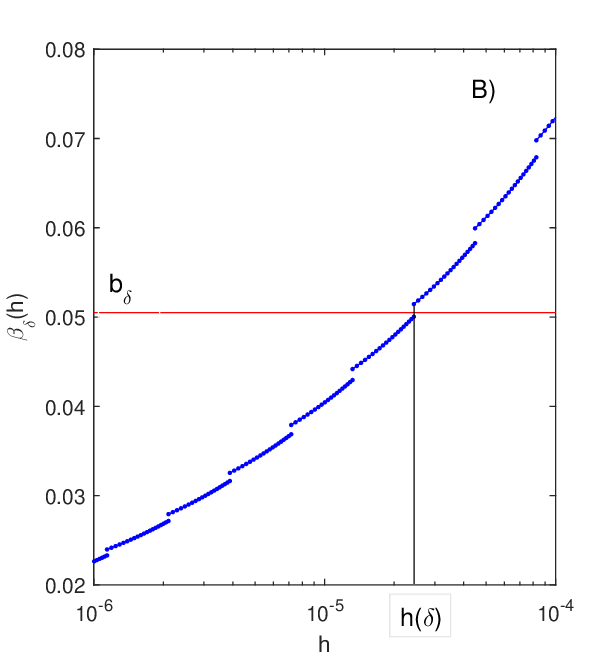}
  \caption{{\small A) To Sect.\ref{sec4}. Matrix singular values \eqref{Poisson}. B) To Sect.\ref{Sec4}. Solution $h(\delta)$ of the equation $\bar\beta_{\delta }(h)=\frac{b_\delta}{\left\Vert \bar{u}\right\Vert _{2}}$ c $b_{\delta }=\sqrt{\mu _{\delta }^{2}+\delta ^{2}}$ for SLAE with matrix \eqref{Poisson}.}}
  \label{fig1}
\end{figure}

In some practical cases, the SLAE matrix is known exactly, often in the form of an analytical expression. In addition, it can be ill-conditioned or even degenerate. Then there are reasons to replace it with a matrix close to it with a better condition number and use the latter instead of the exact one. This can be illustrated by the following example.
\begin{example}\label{ex1}
Application of TSVD method to improve matrix conditioning.
\end{example}
In problems related to potential fields continuation, SLAE with a matrix of the form
\begin{equation}\label{Poisson}
\bar{A}=\left[ \frac{1}{(x_{i}-y_{j})^{2}+H^{2}_0}\right] _{i,j=1}^{m,n},
\end{equation}
must sometimes be solved. Here $\{x_{i}\},\{y_{j}\}$ are uniform grids on the segment $[-1,1]$. For $m=1991,n=2001$ and $H_0=0.1$ the matrix has a condition number $\nu_s=\rho_1/\rho_M\approx 3.37 \cdot 10^{19}$, i.e. is extremely poorly conditioned.
This is due to the specific (exponential) order of decrease of the singular values of the matrix (see Fig.\ref{fig1}A).
By applying the TSVD method to the matrix, i.e., for example, replacing its singular values that satisfy the condition $\rho_k<h$ with zeros, it is possible, with an adequate choice of the value $h$, to improve the condition number. Thus, for $h=10^{-10}$ we obtain a better condition number: $\nu_s\approx 2.48 \cdot 10^{9}$. Accordingly, the stability of numerical solutions of SLAEs with such a matrix improves. However, it is not clear how to constructively find $h$, since the proximity estimate for the matrices, $h$, is unknown in the formulation under consideration. $\square$

We propose the following scheme for solving SLAEs using an approximate matrix with an improved condition number. Let the SVD of the exact matrix be given: $\bar{A}=U\bar{R}V^{T}$. Here $U,V$ are orthogonal matrices of size $m\times m$ and $n\times n$, respectively, and $\bar{R}$ is a diagonal matrix of size $m\times n$ containing the singular values of the matrix $\bar{A}$, sorted in non-increasing order: \[\bar{R}=\mathrm{diag}\left[ \bar{\rho}_{1},\bar{\rho}_{2},...,\bar{\rho}_{\bar r},0,...,0\right] ; \bar{\rho}_{k}\geq \bar{\rho}_{k+1}\geq 0, k=1,2,...,M-1;\,\bar r=\mathrm{Rg}~\bar{A}\leq M.\] We assume that the exact right-hand side $\bar u $ of the SLAE and the matrix $\bar A$ are such that $\|\bar{u}\| \neq 0$, $\|\bar A\|\neq 0$. Then the approximate right-hand side of the SLAE, $u_{\delta }$, satisfies the condition $\|u_{\delta }\| \neq 0$ for sufficiently small values of $\delta$. It is these $\delta$ that we consider further.

\textbf{Scheme} of solving the SLAE with data $\left\{ \bar{A},u_{\delta }\right\} $.

1) For each $h,~0\leq h\leq H=\mathrm{const} $, we look for approximate matrices of the form $\tilde{A}_{h}=U\tilde{R}_{h}V ^{T}$, where \[\tilde{R}_{h}=\mathrm{diag}\left[ \tilde{\rho}_{1}x_{1}(h),\tilde{\rho}_{2}x_{2}(h),...,
\tilde{\rho}_{\bar r}x_{\bar r}(h),0,...,0\right]. \]
The choice of functions $x_{k}(h),\,k=1,...,\bar r$ will be discussed below. The number $H$ must be specified in each method under consideration.

2) We introduce the function $\beta_\delta (h) = \left\| {\bar A\tilde A_h^ + {u_\delta } - {u_\delta }} \right\|$, $0\leq h\leq H$, and solve the equation $\beta^2 _{\delta }(h) =\delta^2 +\beta^2_{\delta }(0)$. We denote its solution as $h(\delta )>0$. Questions about the solvability of the equation will be studied further.

3) We find the matrix $\tilde{A}_{h(\delta )}$ and use it to calculate the approximate solution of the SLAE: $z_{\delta }=\tilde{A}_{h(\delta )}^{+}u_{\delta }$.

Let's make \textbf{assumptions} regarding the functions $x_{k}(h)$: for all $k=1,...,M$ the following requirements must be met.

A) $1<x_{k}(h)\leq c_k=\mathrm{const}$ при $0< h\leq H$;

B) $x_{k}(+0)=x_{k}(0)=1,~x_{k}(H)=0$;

{C) functions $x_{k}(h)$ are left continuous for all $h\in (0,H]$;}

D) the functions $\theta \lbrack x_{k}(h)]$ are non-increasing as $0\leq h\leq H$.

Let us note some simple consequences from these assumptions.

E1) $x_{k}(h)=1~\Longleftrightarrow ~h=0$;
E2) $0\leq \theta \lbrack x_{k}(h)]\leq 1$;
E3) $\underset{h\rightarrow +0}{\lim }\theta \lbrack x_{k}(h)]=\theta \lbrack x_{k}(0)]=1$;
E4) for any $h_{0},~h_{0}\in \lbrack 0,H]$, the equality holds
{\[\overline{\underset{h\rightarrow h_{0}}{\lim }}~\theta \lbrack x_{k}(h)]= \theta \lbrack x_{k}(h_{0})].\]}
We also note the inequality following from E2) and the definition of the function $\theta$
\begin{equation}\label{eq4}
\left\Vert \tilde{A}_{h}^{+}\right\Vert ^{2}=\sum\limits_{k=1}^{\bar r}\bar{\rho}_{k}^{-2}\theta \lbrack x_{k}^{2}(h)]\leq \sum\limits_{k=1}^{\bar r}\bar{\rho}_{k}^{-2}=\left\Vert \bar{A}^{+}\right\Vert ^{2}.
\end{equation}

\begin{theorem}\label{thm1}
Let conditions A) - D) be satisfied and, in addition, for each $\delta ,~0<\delta <\delta _{0}=\mathrm{const}$, the inequality $\left\Vert u_{\delta }\right\Vert >\mu _{\delta }=\left\Vert \bar{A}\bar{A}^{+}u_{\delta }-u_{\delta }\right\Vert $ is valid. Then:

1) the function $\beta _{\delta }(h)$ non-decreasing as $h\in \lbrack 0,H]$ and is left continuous at every point $h>0$;

2) $\beta^2 _{\delta }(+0)=\mu _{\delta }^{2},~\beta^2_{\delta }(\left\Vert \bar{A}\right\Vert )=\left\Vert u_{\delta }\right\Vert ^{2}$;

3) the equation $\beta^2_{\delta } (h) = \delta^2+\mu^2_\delta$ has a generalized solution $h(\delta )>0$;

4) $h(\delta )\rightarrow 0$ при $\delta \rightarrow 0$.
\end{theorem}
\begin{remark}\label{rem1}
In item 3) the equation with the monotonic function $\beta_{\delta } (h)$ is solved. Its generalized solution $h(\delta )$ is a point for which the inequalities hold
\[\underset{h\rightarrow h(\delta )}{\underline{\lim }}\beta^2 _{\delta }(h)=\beta^2_{\delta }(h(\delta )-0)\leq \mu _{\delta }^{2}+\delta ^{2}~\leq \beta^2_{\delta }(h(\delta )+0)=\underset{h\rightarrow h(\delta )}{\overline{\lim }}\ \beta ^2_{\delta }(h).\]
Such equations arose in the study of the discrepancy principle for solving ill-posed inverse problems (see, for example, \cite{tikh2}). In what follows we will use the part of this inequality that takes into account the left continuity of the function $\beta_{\delta } (h)$:
\begin{equation}\label{sol_bet}
\beta ^2_{\delta }(h(\delta ))=\beta ^2_{\delta }(h(\delta )-0)\le \delta^2+\mu^2_\delta.
\end{equation}
\end{remark}
Proof. First, we prove 1). Equality
\[\beta ^2_{\delta }(h)=\left\Vert \bar{A}^{+}\tilde{A}_{h}^{+}u_{\delta }-u_{\delta }\right\Vert ^{2}=\left\Vert U\bar{R}V^{T}V\tilde{R}_{h}^{+}U^{T}u_{\delta }-UU^{T}u_{\delta }\right\Vert ^{2}=\left\Vert \bar{R}\tilde{R}_{h}^{+}v_{\delta }-v_{\delta }\right\Vert^2 ,\]где $v_{\delta }=U^{T}u_{\delta }=(v_{k})_{k=1}^{m}$ fulfilled.
Next, taking into account the form of the matrices $\bar{R}$ and $\tilde{R}_{h}^{+}$, we obtain
\begin{gather}
\beta ^2_{\delta }(h) =\sum\limits_{k=1}^{m}(\bar{\rho}_{k}\theta \lbrack \bar{\rho}_{k}x_{k}(h)]-1)^{2}v_{k}^{2}=\sum\limits_{k=1}^{\bar r}
(\bar{\rho}_{k}\bar{\rho}_{k}^{-1}\theta \lbrack x_{k}(h)]-1)^{2}v_{k}^{2}+\sum\limits_{k=\bar r+1}^{m}v_{k}^{2}=\nonumber\\
=\sum\limits_{k=1}^{\bar r}(1-\theta \lbrack x_{k}(h)])^{2}v_{k}^{2}+\mu _{\delta }^{2}.\label{mull}
\end{gather}
From this, as well as from properties D) and E2), the monotonicity (non-decreasing) of the function $\beta _{\delta }(h)$ follows.

The equalities indicated in part 2) of the theorem follow from properties B), E3) and E4). The existence of a generalized solution to the equation from item 3) follows from the monotonicity of the function $\beta _{\delta }(h)$, the proven relations of item 2) and the condition $\left\Vert u_{\delta }\right\Vert >\mu _{\delta }$ of theorem (for a similar statement, see \cite{leonov1}).

Let us prove the convergence from part 4). Suppose that $h(\delta)\nrightarrow 0$ for $\delta \rightarrow 0$. Then there is a sequence $\{\delta _{N}\},\delta _{N}\rightarrow 0$ for which $h(\delta _{N})=h_{N}\rightarrow h_{0}> 0$ as $N\rightarrow \infty$. From the equality \eqref{mull} and the property \eqref{sol_bet} of solutions $h(\delta _{N})=h_{N}$ it follows:
\[\beta^2 _{\delta }(h_{N})=\sum\limits_{k=1}^{\bar r}(\theta \lbrack x_{k}(h_{N})]-1)^{2}(v_{\delta _{N}})_{k}^{2}+\mu _{\delta _{N}}^{2}\leq \delta _{N}^{2}+\mu _{\delta _{N}}^{2}\Longrightarrow \sum\limits_{k=1}^{\bar r}(\theta \lbrack x_{k}(h_{N})]-1)^{2}(v_{\delta _{N}})_{k}^{2}\leq \delta _{N}^{2}.\]
Passing to the limit at $N\rightarrow \infty $ in the last inequality and taking into account E4), we obtain: \[\sum\limits_{k=1}^{\bar r}(\theta \lbrack x_{k}(h_{0})]-1)^{2}\bar{v}_{k}^{2}=0.\]
where $(\bar{v}_{k})_{k=1}^{m}=U^{T}\bar{u}$. It follows that there is a number $k$ for which $\theta \lbrack x_{k}(h_{0})]=1$, because $\left\Vert \bar{u}\right\Vert =\left\Vert U^{T}\bar{u}\right\Vert \neq 0$. But then, by properties D), E1), E2) it turns out that $x_{k}(h_{0})=1\Longrightarrow h_{0}=0$, i.e. $h_{N}\rightarrow 0$. The resulting contradiction proves the theorem.

\begin{theorem}\label{thm2}
The approximate solution of the SLAE, $z_{\delta }=\tilde{A}_{h(\delta )}^{+}u_{\delta }$, converges to the normal pseudo-solution $\bar z$ for $\delta \rightarrow 0$ .
\end{theorem}
Proof. From the form of the matrices $\tilde{A}_{h(\delta )},\bar{A}$ it follows that \[\left\Vert \tilde{A}_{h(\delta )}-\bar{A }\right\Vert ^{2}=\sum\limits_{k=1}^{\bar r}\bar{\rho}_{k}^{2}\left( x_{k}(h(\delta ))-1\right) ^{2},\] and by virtue of property B) and the proven convergence $h(\delta )\rightarrow 0$ the equality $\underset{\delta \rightarrow 0}{\lim }x_{ k}(h(\delta ))=x_{k}(0)=1$ is true. Therefore $\tilde{A}_{h(\delta )}\rightarrow \bar{A}$. Now applying corollary \ref{cor2} taking into account \eqref{eq4}, we obtain $\tilde{A}_{h(\delta )}^{+}\rightarrow \bar{A}^{+}$ and $z_ {\delta }=\tilde{A}_{h(\delta )}^{+}u_{\delta }\rightarrow \bar{z}=\bar{A}^{+}\bar{u}$ .$\square$

Now we study the question of improving the condition number of the original matrix $\bar{A}$ when passing to the matrix $\tilde{A}_{h(\delta )}$.
\begin{theorem}\label{thm3}
Let $r(\delta )$ be the rank of the matrix $\tilde{A}_{h(\delta )}$. Suppose that for each $k$, $1\le k\le \bar r$, the relation $x_{k}(h)\sim 1+a_{k}h$ holds when $h\rightarrow +0$ , with $ a_{1}<a_{2}<...<a_{\bar{r}}$. Then for $\delta \rightarrow 0$ the following estimate is valid: $\nu _{s}(\tilde{A}_{h(\delta )})\sim \frac{\bar{\rho}_{1}}{\bar{\rho}_{r(\delta )}}(1-h(\delta )(a_{r(\delta )}-a_{1}))<\nu _{s}(\bar{A})$.
\end{theorem}
Proof. From the form of singular values of the matrix $\tilde{A}_{h(\delta )}$, the equality $\nu _{s}(\tilde{A}_{h(\delta })=\frac{\bar){\rho}_{1}x_{1}(h(\delta ))}{\bar{\rho}_{r(\delta )}x_{r(\delta )}(h(\delta )) }$ follows, and $r(\delta )\leq \bar{r}$. Hence, taking into account the conditions of the theorem and the convergence $h(\delta )\rightarrow 0$, we obtain a chain of relations
\begin{multline*}
\nu _{s}(\tilde{A}_{h(\delta )}) =\frac{\bar{\rho}_{1}x_{1}(h(\delta ))}{\bar{\rho}_{r(\delta )}x_{r(\delta )}(h(\delta ))}\sim \frac{\bar{\rho}_{1}}{\bar{\rho}_{r(\delta )}}\frac{1+a_{1}h(\delta )}{1+a_{r(\delta )}h(\delta )}\sim \frac{\bar{\rho}_{1}}{\bar{\rho}_{r(\delta )}}(1+a_{1}h(\delta ))(1-a_{r(\delta )}h(\delta ))\sim \\
\sim \frac{\bar{\rho}_{1}}{\bar{\rho}_{r(\delta )}}(1-h(\delta )(a_{r(\delta )}-a_{1}))<\frac{\bar{\rho}_{1}}{\bar{\rho}_{r(\delta )}}\leq \frac{\bar{\rho}_{1}}{\bar{\rho}_{\bar{r}}}=\nu_{s}(\bar A),
\end{multline*}
that proves the theorem.

Thus, the condition number of the matrix of the proposed method is better than the condition number of the original exact matrix, at least for sufficiently small $\delta$. If, in addition to the conditions of the theorem, the exact matrix has a significant number of singular values close to zero, as in the example \ref{ex1}, then it may turn out that $\bar{\rho}_{r(\delta )}>>\bar {\rho}_{\bar{r}}$, and then $\nu _{s}(\tilde{A}_{h(\delta )})<<\nu _{s}(\bar{ A})$. This will be illustrated below.

\section{A specific method that implements the proposed scheme for solving SLAE}\label{Sec4}
For simplicity, we assume that the singular values of the matrix $\bar A$ are different: $\bar{\rho}_{k}>\bar{\rho}_{k+1},~k=1,2,...,M-1$.

As a central example, consider {method using minimal pseudoinverse matrices with condition number improvement} or briefly, \textbf{MPMI method}. This method uses the quantities \[x_{k}(h)=\left\{ \bar{x}_{k}(h),~0\leq h\leq h_{k};~0,~h> h_{k}\right\}\] for $1\le k\le \bar{r}$, where $\bar{x}_{k}(h)$ is the solution to the equation $\bar{x}_{k}^{4}-\bar{x}_{k}^{3}=h\bar{\rho}_{k}^{-4},~\bar{x}_{k}\in \lbrack 1,\frac{3}{2}]$, а $h_{k}=\frac{27}{16}\bar{\rho}_{k}^{4}$. Тогда \[\theta \left[ x_{k}(h)\right] =\left\{ \frac{1}{\bar{x}_{k}(h)},~0\leq h\leq h_{k};~0,~h>h_{k}\right\}. \] Let us also define $H$: $H>h_1=\frac{27}{16}\bar{\rho}_{1}^{4}$. From the definition of the function $x_{k}(h)$ and from the form of the function $\bar{x}_{k}(h)$ the following properties are obtained: A) $1\leq x_{k}(h)\leq \frac{3}{2}$, $1\le h\le H$; B) $x_{k}(0)=x_{k}(+0)=\bar{x}_{k}(0)=1$; C) $x_{k}(h)$ is left continuous for $0<h\le H$. Property D) is also true. It follows from the increase in the function $\bar{x}_{k}(h)$ as $h\in \lbrack 0,h_{k}]$ (see \cite{leonov1,tikh2}). Then $\theta \left[ x_{k}(h)\right] =\frac{1}{\bar{x}_{k}(h)}>0$ decreases as $h\in \lbrack 0, h_{k}]$ and is equal to zero for $h>h_{k}$. Analyzing the asymptotic behavior of the functions $x_{k}(h)$ for $h\rightarrow +0$, we can verify that $x_{k}(h)\sim 1+h\bar{\rho}_{k}^{-4}$. This means that the conditions of the theorem \ref{thm3} are satisfied with $a_{k}=\bar{\rho}_{k}^{-4}$, and $a_{k}=\bar{\rho}_{k }^{-4}<\bar{\rho}_{k+1}^{-4}=a_{k+1}$ for all admissible $k$.

Thus, theorems \ref{thm2} and \ref{thm3} are valid, guaranteeing the convergence of the MPMI method and the improvement of the condition number in it. In some cases, this number differs significantly from $\nu _{s}(\bar{A})$. For example, the theoretically possible case is $x_{r(\delta )}(h(\delta ))=x_{r(\delta )}(h_{r(\delta )})=\frac{3}{2}$ . It is realized when the generalized solution of the equation $\beta _{\delta }^{2}(h)=\delta ^{2}+\mu _{\delta }^{2}$ is the discontinuity point $h_{r( \delta )}$ functions $\beta _{\delta }(h)$. Then, as in theorem \ref{thm3}, we obtain the relations
\begin{equation}\label{best}
\nu _{s}(\tilde{A}_{h(\delta )})=\frac{2}{3}\frac{\bar{\rho}_{1}x_{1}(h(\delta ))}{\bar{\rho}_{r(\delta )}}\sim \frac{2}{3}\frac{\bar{\rho}_{1}}{\bar{\rho}_{r(\delta )}}(1-h(\delta )(\bar{\rho}_{r(\delta )}^{-4}-\bar{\rho}_{k}^{-4}))<\frac{2}{3}\frac{\bar{\rho}_{1}}{\bar{\rho}_{r(\delta )}}\leq \frac{2}{3}\nu _{s}(\bar{A}),
\end{equation}
and in this case, the condition number of the matrix used to solve the SLAE is improved by at least one and a half times.

Next, we compare the MPMI method with some others.

\textbf{TSVD method.} Here, instead of $\bar{A}$, the matrix $\tilde{A}_{r(\delta )}=U\tilde{R}_{r(\delta )}V^{T}$ is used, in which $ \tilde{R}_{r(\delta )}=\mathrm{diag}\left[ \bar{\rho}_{1},\bar{\rho}_{2},...,\bar {\rho}_{r(\delta )},0,...,0\right] ~$. Its rank $r(\delta )$ is found as a solution to the equation $\beta^2 _{\delta }(r)=\mu _{\delta }^{2}+\delta ^{2}$ with monotonically non-decreasing function $\beta^2_{\delta }(r)=\sum\limits_{k=r+1}^{m}v_{k}^{2}$. We can verify that $r(\delta )=r$ for sufficiently small $\delta $. Thus, the inequality $\nu _{s}(\tilde{A}_{r(\delta )})=\frac{\bar{\rho}_{1}}{\bar{\rho}_{ r(\delta )}}\leq \frac{\bar{\rho}_{1}}{\bar{\rho}_{r}}=\nu _{s}(\bar{A})$ becomes equality for small $\delta $, and the TSVD method does not improve the condition number.

\textbf{Variants of the regularization method. }In the Tikhonov regularization (TR), an approximate solution of the SLAE is sought in the form $z_{\delta }=\left( \alpha (\delta )I+\bar{A}^{\ast }\bar{A}\right) ^{- 1}\bar{A}^{\ast }u_{\delta }$, where the parameter $\alpha (\delta )>0$ is selected by one of the known methods (see, for example, \cite{tikh1,tikh2,tikh3 ,moroz}). In any case, $\alpha (\delta )\rightarrow 0$ for $\delta \rightarrow 0$. Using the SVD for the matrix $\bar{A}$ we find that $z_{\delta }=VT_{\delta }^{-1}U^{T}u_{\delta }$, where
\[T_{\delta }=\mathrm{diag}\left[ \frac{\alpha (\delta )+\bar{\rho}_{1}^{2}}{\bar{\rho}_{1}},\frac{\alpha (\delta )+\bar{\rho}_{2}^{2}}{\bar{\rho}_{2}},...,\frac{\alpha (\delta )+\bar{\rho}_{r}^{2}}{\bar{\rho}_{r}},0,...,0\right]. \] Then, as in the proof of theorem \ref{thm3}, for $\alpha (\delta )\rightarrow 0$ we obtain
\begin{multline*}
 \nu _{s}(T_{\delta })=\frac{\bar{\rho}_{r}(\alpha (\delta )+\bar{\rho}_{1}^{2})}{\bar{\rho}_{1}(\alpha (\delta )+\bar{\rho}_{r}^{2})}=\frac{\bar{\rho}_{1}}{\bar{\rho}_{r}}\frac{(1+\alpha (\delta )/\bar{\rho}_{1}^{2})}{(1+\alpha (\delta )/\bar{\rho}_{r}^{2})}\sim \nu _{s}(\bar{A})\left( 1+\frac{\alpha (\delta )}{\bar{\rho}_{1}^{2}}\right) \left( 1-\frac{\alpha (\delta )}{\bar{\rho}_{r}^{2}}\right) \sim\\
\sim \nu _{s}(\bar{A})\left( 1-\alpha (\delta )\left( \frac{1}{\bar{\rho}_{r}^{2}}-\frac{1}{\bar{\rho}_{1}^{2}}\right) \right) <\nu _{s}(\bar{A}),
\end{multline*}
taking into account that $\rho _{1}>\rho _{r}$.

Similar calculations cab be carried out for the method of solving SLAEs from the work \cite{moroz}:
\[z_{\delta }=\left( \alpha (\delta )I+\bar{A}^{T}\bar{A}\right) ^{-1}\bar{A}^{T}\bar{A}\bar{A}^{T}\left( \alpha (\delta )I+\bar{A}\bar{A}^{T}\right) ^{-1}u_{\delta }=VM_{\delta }^{-1}U^{T}.\]As a result, we get a similar estimate:
\[\nu _{s}(M_{\delta })\sim \nu _{s}(\bar{A})\left( 1-\alpha (\delta )\left( \frac{1}{\bar{\rho}_{r}^{2}}-\frac{1}{\bar{\rho}_{1}^{2}}\right) \right) ^{2}<\nu _{s}(\bar{A}).\]Thus, these variants of the regularization method theoretically provide, in comparison with the original matrix $\bar A$, an improvement in the condition number for the matrices $T_{\delta }$ and $M_{\delta }$, with the help of which SLAEs are solved. However, these matrices have the same rank as $\mathrm{rank}\bar{A}$, while in the MPMI method the rank $r(\delta )$ of the matrix $\tilde{A}_{h(\delta )}$ can be significantly less than $\bar r$. This possibility is demonstrated below in a numerical experiment. As a consequence, the main terms of the condition numbers $\nu _{s}(T_{\delta })$ and $\nu _{s}(M_{\delta })$, expressed by the term $\frac{\bar{\rho }_{1}}{\bar{\rho}_{\bar r}}$ cannot be better than the leading term of the condition number $\nu _{s}(\tilde{A}_{h(\delta )} )\sim \frac{\bar{\rho}_{1}}{\bar{\rho}_{r(\delta )}}$. In addition, the estimate \eqref{best} for the number $\nu _{s}(\tilde{A}_{h(\delta )})$ may be fair. This guarantees a reduction of the condition number by one and a half times. This feature is also demonstrated below.

\section{Numerical experiments}
Let us illustrate the main features of the proposed scheme for solving SLAE using the example of the MPMI method discussed above. We solve a model SLAE with an ill-conditioned matrix $\bar{A}$ of the form \eqref{Poisson} and an exact solution $\bar z=(1-x^2)\sin(4 \pi x), \,-1 \le x\le 1$. The exact right-hand side of this SLAE is calculated as $\bar{u}=\bar{A}\bar{z}$ and is perturbed by a normally distributed random noise with zero mean so that the approximate right-hand side $u_{\delta }$ satisfies the inequality $\left\Vert u_{\delta }-\bar{u}\right\Vert _{2}\leq \delta \left\Vert \bar{u}\right\Vert _{2}$. The number $\delta $ specified in this way allows one to estimate the error level as a percentage. To solve the considered SLAE with different error levels $\delta $, the MPMI method is used. For comparison, the TSVD and TR methods were also used with the choice of the regularization parameter $\alpha(\delta)$ according to the discrepancy principle \cite{moroz}.

Table 1 shows the results of calculations, namely the accuracy of the obtained approximate solutions $z_{\delta }$: $\Delta =\frac{\left\Vert z_{\delta }-\bar{z}\right\Vert _{2}}{\left\Vert \bar{z}\right\Vert _{2}}$ and condition numbers $\nu _{MPMI}=\nu _{s}(\tilde{A}_{h(\delta )})$, $\nu _{TSVD}=\nu _{s}(\tilde{A}_{r(\delta )})$ and $\nu _{TR}=\nu _{s}(T_{\delta })$.
\begin{center}
Table 1. Comparison of solution accuracies and condition numbers\\ matrices of the MPMI, TSVD and TR methods
\begin{tabular}{|c|c|c|c|c|c|c|}
\hline
$\delta $ & 0.005 & 0.01 & 0.05 & 0.1 & 0.2 & 0.3 \\ \hline
$\Delta _{MPMI}$ & 0.0024 & 0.0043 & 0.0117 & 0.0154 & 0.0333 & 0.0406 \\
\hline
$\Delta _{TSVD}$ & 0.0027 & 0.0052 & 0.0131 & 0.0184 & 0.0346 & 0.0496 \\
\hline
$\Delta _{TR}$ & 0.0082 & 0.0108 & 0.0269 & 0.0358 & 0.0495 & 0.0989 \\
\hline
$\nu _{MPMI}$ & 20.972 & 20.971 & 10.353 & 10.353 & 10.353 & 5.6134 \\ \hline
$\nu _{TSVD}$ & 33.421 & 33.420 & 15.530 & 15.530 & 15.530 & 8.4172 \\ \hline
$\nu _{TR}$ & 2.3$\cdot 10^{12}$ & 5.9$\cdot 10^{12}$ & 3.3$\cdot 10^{13}$
& 7.7$\cdot 10^{13}$ & 2.6$\cdot 10^{14}$ & 5.6$\cdot 10^{14}$ \\ \hline
\end{tabular}
\end{center}
\vspace{1mm}
The table shows that the MPMI method has the best accuracy for all considered levels $\delta$ of data disturbance. We especially note the high accuracy of this method for large disturbances: $\delta =0.1,0.2,0.3$. Also, the table clearly shows how much the condition numbers of the matrices of the MPMI method are smaller than the condition numbers of other methods and, especially, the condition numbers of the original SLAE matrix, $\nu_s\approx 3.37 \cdot 10^{19}$. Note also that for the method from the work \cite{moroz} the results are very close to their analogues from the TR method and therefore are not included in the table.

Calculations also showed that a significant decrease in the condition number in the MPMI method (see \eqref{best}) occurs quite often.
For example, for $\delta =0.05$ it turns out that $h(\delta )$ is a discontinuity point of the function $\beta _{\delta }(h)$, and therefore for such $h(\delta )$ the estimate \eqref{best} is true (see fig.\ref{fig1}B with $\bar{\beta}_{\delta }(h)=\beta_{\delta }(h)/\|\bar u\|_2$). Effects of this kind also occurred for other $\delta $.

In conclusion, we note that similar results were also obtained in numerical experiments with other SLAEs that have ill-conditioned matrices.

This work was supported by the Russian Science Foundation (grant no. 23-41-00002).


\begin{thebibliography}{2}

\bibitem{fox} Fox L. An introduction to numerical linear algebra, Oxford Univ. Press, NY, 1965.

\bibitem{fors} Forsythe G.E. and Moler C. Computer Solution of Linear Algebraic Systems, Prentice Hall,
Englewood Cliffs, NJ, 1967.

\bibitem{lowson} Lawson C.L. and Hanson R.J. Solving Least Squares Problems, Prentice Hall,
Englewood Cliffs, NJ, 1974.

\bibitem{golub1} Golub G.H. and Van Loan C.F. Matrix Computations, The John Hopkins University Press, 1989

\bibitem{voev} Voevodin Val.V. Numerical methods, Nauka, Moscow, 1977 (in Russian)
.
\bibitem{tyrt} Tyrtyshnikov E.E. Methods of numerical analysis, Publishing center Academy, Moscow, 2007 (in Russian).

\bibitem{tikh1} Tikhonov A.N. and Arsenin V.Y. Solution of Ill-Posed Problems, Wiley, New York,
1977.

\bibitem{TSVD1} Varah J.M. On the numerical solution of ill-conditioned linear systems with applications to ill-posed problems,
SIAM J. Num. Anal., V.10, 257–267 (1973).

\bibitem{TSVD2} Engl H.W., Hanke M. and Neubauer A. Regularization of Inverse Problems,
Kluwer, Dordrecht, 1996.

\bibitem{golub2} Golub G.H. Least squares, singular values and matrix approximation,
Aplikace Matematiky, V.13, 44–51 (1968).

\bibitem{leonov1}  Leonov A.S. The minimal pseudo-inverse matrix method,
U.S.S.R. Comput. Math. Math. Phys., V.27, No.4, 107–117 (1987).

\bibitem{tikh2} Tikhonov A.N., Leonov A.S. and Yagola A.G. Nonlinear Ill-Posed Problems, V.2,
 Chapman and Hall, London, 1998.

\bibitem{leonov_book} Leonov A.S. Solving Ill-Posed Inverse Problems: Outline of Theory,
Practical Algorithms and Demonstrations in MATLAB, Librokom, Moscow, 2010 (in Russian).

\bibitem{alb} Albert A. Regression and the Moor-Penrose pseudoinverse, Academic Press,
New York and London, 1972.

\bibitem{tikh3} Tikhonov A.N., Goncharsky A.V., Stepanov V.V. and Yagola A.G. Numerical Methods for the Solution of Ill-Posed Problems, Mathematics and Its Applications,  Vol. 328, Springer, Dordrecht, 1995.

\bibitem{moroz}Morozov V.A.  Methods for Solving Incorrectly Posed Problems, Springer-Verlag, New York, 1984.






\end{thebibliography}
\end{document}